\documentclass[a4paper,10pt]{article}
\usepackage[utf8x]{inputenc}
\usepackage{amsmath,amsthm}
\usepackage{amsfonts}
\usepackage{amssymb}
 \usepackage{mathtools}
 \usepackage[all]{xy}
 \usepackage{tikz}
 \usetikzlibrary{positioning}
 \usetikzlibrary{decorations.pathreplacing}
 \usepackage{textcomp}

\theoremstyle{definition}
 \newtheorem{thm}{Theorem}
 \newtheorem*{thm*}{Theorem}
 \newtheorem{cor}[thm]{Corollary}
 \newtheorem{lem}[thm]{Lemma}

 \theoremstyle{definition}
 
 \theoremstyle{definition}
 \newtheorem*{defn*}{Definition}

\theoremstyle{remark}
\newtheorem{rmk}{Remark}

\theoremstyle{definition}

\newtheoremstyle{named}{}{}{\itshape}{}{\bfseries}{.}{.5em}{\thmnote{#3}}
\theoremstyle{named}
\newtheorem*{namedtheorem}{Theorem}

\newtheoremstyle{namedproof}{}{}{}{}{\bfseries}{.}{.5em}{\thmnote{#3}}
\theoremstyle{namedproof}
\newtheorem*{namedproof}{namedproof}

\theoremstyle{definition}
\newtheorem*{conj}{Conjecture}

\newcommand{\bbb}[1]{\textbf{#1}}
\newcommand{\vvert}{\ \vert\ }
\newcommand{\norm}[1]{{\vert #1 \vert}}

\newcommand{\wt}[1]{\widetilde{#1}}
\newcommand{\wh}[1]{\widehat{#1}}
\newcommand{\cwr}{\mbox{\small\textcircled{$\wr$}}}

\newcommand{\ul}[1]{\boldsymbol{#1}}
\newcommand{\defeq}{\vcentcolon=}

\newcommand{\comm}[1]{}

\newcommand*{\transp}[2][-3mu]{\ensuremath{\mskip1mu\prescript{\smash{\mathrm t\mkern#1}}{}{\mathstrut#2}}}
\newcommand{\mS}{\mathcal{S}}

\renewcommand{\d}{\mathrm{d}}
\renewcommand{\S}{\mathrm{Sym}}
\renewcommand{\a}{\alpha}
\renewcommand{\b}{\beta}

\title{Finite generation of  \\ iterated wreath products in product action}
\author{Matteo Vannacci}

\begin{document}

\maketitle

\begin{abstract}
 Let $\mathcal{S}$ be a sequence of finite perfect transitive permutation groups with uniformly bounded number of generators. We prove that the infinitely iterated wreath product in product action of the groups in $\mS$ is topologically finitely generated, provided that the actions of the groups in $\mS$ are not regular. We prove that our bound has the right asymptotic behaviour. We also deduce that other infinitely iterated mixed wreath products of groups in $\mS$ are finitely generated. Finally we apply our methods to find explicitly two generators of infinitely iterated wreath products in product action of special sequences $\mS$.
\end{abstract}


\phantom{a}

\noindent\bbb{Notation 1.} All the actions will be right actions and $e$ will stand for the identity element of a group. We will write $\ul{n}=\{1,\ldots,n\}$. 

\phantom{a}

\noindent\bbb{Notation 2.} Let $d$ be an integer. In all of this paper we will denote by $\mS$ a sequence of finite transitive permutation groups $\{S_k \le \S(m_k)\}_{k\in \mathbb{N}}$ such that each $S_k$ is perfect and at most $d$-generated as an abstract group.

\phantom{a}

\section{Introduction}
Infinitely iterated wreath products have been widely studied in the past (see \cite{bhattacharjee,MR2727305,quick:probabilisticgeneration,woryna:cyclic,woryna:abelian}) and their generation properties proved to be of great interest.
In particular, \cite[Theorem 1]{MR2727305} states that an infinitely iterated permutational wreath product of finite $d$-generated transitive permutation groups is finitely generated if and only if its abelianization is finitely generated. In this paper we prove two parallel results. We prove that an \emph{infinitely iterated exponentiation} and an \emph{infinitely iterated mixed wreath product of stride at most $m$} of finite $d$-generated perfect transitive permutation groups are topologically finitely generated under certain conditions.

\begin{defn*}
The abstract wreath product $A\wr B$ of two permutation groups $A\le \S(\ul{m})$ and $B\le \S(\ul{n})$ can be considered as a permutation group in two (in general non-equivalent) ways. One way is to consider $A\wr B$ acting on $\ul{m}^n$ with the product action, we will write $A \cwr B\le \S(\ul{m}^n)$ for this permutation group and we will call it the \bbb{exponentiation} of $A$ by $B$. The second way is to consider $A\wr B$ acting on $\ul{m}\times \ul{n}$ with the permutational wreath action, we will write $A \wr B\le \S(\ul{m}\times\ul{n})$ for this permutation group and we will call it the \bbb{permutational wreath product} of $A$ by $B$.
\end{defn*}
 For other properties of the exponentiation and permutational wreath product we refer to \cite{MR1409812} and references therein.
\begin{defn*}
The \bbb{iterated exponentiation} $\wt{S}_n \le \S(\ul{\wt{m}}_n)$ of the groups in the sequence $\mathcal{S}$ is the permutation group inductively defined by: $\wt{m}_1 = m_1$, $\wt{S}_1= S_1 \le \S(\wt{\ul{m}}_1)$ and $\wt{m}_{k} = m_k^{\wt{m}_{k-1}}$, $\wt{S}_{k} = S_{k} \cwr \wt{S}_{k-1} \le \S(\ul{\wt{m}}_{k})$ for $k\ge 2$. The groups $\wt{S}_k$, together with the projections $\wt{S}_{k} \rightarrow \wt{S}_{k-1}$, form an inverse system of finite groups. We will call the profinite group $\varprojlim \wt{S}_k$ the \bbb{infinitely iterated exponentiation} of the groups in $\mS$.
\end{defn*}

In Section \ref{sec:A} we are going to prove our main result.

\begin{namedtheorem}[Theorem A]
 Let $d$ be an integer. Let $\mS = \{S_k \le \S(m_k)\}_{k\in \mathbb{N}}$ be a sequence of finite transitive permutation groups such that each $S_k$ is perfect and at most $d$-generated as an abstract group. Suppose that for every $k \in \mathbb{N}$ there exist elements $i,j\in \ul{m}_k$ such that $\mathrm{St}_{S_k}(i) \neq \mathrm{St}_{S_k}(j)$. Then the infinitely iterated exponentiation of the groups in $\mS$ is topologically finitely generated. In particular, we exhibit a set of $d + \mathrm{d}(S_1)$ generators for $\varprojlim \wt{S}_k$.
\end{namedtheorem}

Using the same methods we can improve our bound for a sequence $\mS$ of perfect 2-generated perfect groups (see Corollary \ref{cor:3gen}). 

In Section \ref{sec:B} we are going to use \cite[Theorem 3.1]{kaluzhnin} to exchange the exponentiation and the permutational wreath product to obtain the finite generation of iterated wreath products with ``mixed'' action.

\begin{defn*}
Let $\{k_n\}_{n\in \mathbb{N}}$ be an increasing sequence of integers.  Define the new sequence $\{G_n\le \S(\ul{r}_n)\}_{n\in \mathbb{N}}$ of finite perfect transitive permutation groups starting from the groups in $\mS$ in the following way: $G_1 = S_1$, $G_{k+1} = S_{k+1}\wr G_k$ for $k\in \mathbb{N}\smallsetminus \{k_n\}_{n\in \mathbb{N}}$ and $G_{k_n} = S_{k_n}\cwr G_{k_n-1}$ for every $n \in \mathbb{N}$. The permutation groups $G_n$ are called \bbb{iterated mixed wreath product of type} $(\mS,\{k_n\}_{n\in \mathbb{N}})$. 
 
 Let $m$ be an integer, if the sequence $\{k_n\}_{n\in \mathbb{N}}$ is such that $k_{n+1}-k_n \le m$ for every $n\in \mathbb{N}$, we say that the iterated mixed wreath product $G_n$ of type $(\mS,\{k_n\}_{n\in \mathbb{N}})$ has \bbb{stride at most $m$}.
 
 The groups $G_n$, together with the projections $G_n \rightarrow G_{n-1}$, form an inverse system of finite groups. We say that the profinite group $\varprojlim G_n$ is an \bbb{infinitely iterated mixed wreath product of type $(\mS,\{k_n\}_{n\in \mathbb{N}})$}. If the groups $G_n$ have stride at most $m$ we say that $\varprojlim G_n$ has stride at most $m$.
\end{defn*}

An infinitely iterated exponentiation is an infinitely iterated mixed wreath product of stride at most one. Our second main result is the following.

\begin{namedtheorem}[Theorem B]
 Let $d$ be an integer. Let $\mS = \{S_k \le \S(m_k)\}_{k\in \mathbb{N}}$ be a sequence of finite transitive permutation groups such that each $S_k$ is perfect and at most $d$-generated as an abstract group. Suppose that for every $k \in \mathbb{N}$ there exist elements $i,j\in \ul{m}_k$ such that $\mathrm{St}_{S_k}(i) \neq \mathrm{St}_{S_k}(j)$. Let $G = \varprojlim G_n$ be an infinitely iterated mixed wreath product of type $(\mS,\{k_n\}_{n\in \mathbb{N}})$ of stride at most $m$. Then $G$ is topologically finitely generated. In particular the profinite group $G$ is generated by at most $2md$ elements. 
\end{namedtheorem}
The hypotheses of Theorem B can be weakened in two ways (see Remark \ref{rmk:weak}).

We conclude with Section \ref{sec:special} where we use the techniques of this paper to find the minimal number of generators of infinitely iterated exponentiations of particular sequences $\mS$ (see Corollary \ref{cor:special}).

\section{Proof of Theorem A}\label{sec:A}
First we find a lower bound for the minimal number of generators of a wreath product of perfect non-simple groups showing that the bound given by Theorem A could possibly be improved only by a multiplicative and an additive constant. We will denote by $\transp{x}$ the transpose of a vector $x\in \mathbb{Z}^n$.

\begin{lem}\label{lem:lowerbound}
Let $N$ be a natural number. Let $A$ be a finite simple group and let $B\le \mathrm{Sym}(n)$ be a finite perfect permutation group. Then
\[
  \d(A^N\wr B) \ge \max\left\{\frac{1}{n}\left( \mathrm{d}\left(A^N\right) - \d(A)-1\right),\d(B)\right\}.
\]
\begin{proof}
Set $G = A^N\wr B = (A^{N})^{n} \rtimes B$ and $\d(G) =d$. It is clear that $d\ge \d(B)$, since $B$ is a quotient of $G$. Let $$\gamma_j = ((x_{11}^{j},\ldots,x_{1N}^j),\cdots,(x_{n1}^j,\ldots,x_{nN}^j)) \sigma_j = (\gamma_{1}^{(j)},\cdots, \gamma_{n}^{(j)}) \sigma_j\in A^{N}\wr B,$$ for $j=1,\ldots,d$, be generators for $G$. Form the $N \times nd$ matrix $M$ with entries $M_{l,n(i-1) +j} = x_{jl}^i$ for $i=1,\ldots,d$, $j=1,\ldots,n$ and $l=1,\ldots,N$. For every number $m\in\{1,\ldots,nd\}$ there exist unique $i\in\{1,\ldots,d\}$ and $j\in \{1,\ldots,n\}$ such that $m=n(i-1)+j$ and the $[n(i-1)+j]$-th column of $M$ is the vector $\transp{\gamma^{(i)}_j}$.
\[
  M=\left(\transp{\gamma_1^{(1)}},\ldots,\transp{\gamma_n^{(1)}},\ldots\ldots,\transp{\gamma_1^{(d)}},\ldots,\transp{\gamma_n^{(d)}}\right).
\]
Our goal is to show that $N\le \norm{A}^{nd}$. Suppose by contradiction that $N>\norm{A}^{nd}$. Then, since the $x_{jl}^i$ are elements of $A$, we would have that two rows of $M$ are equal. Without loss of generality we can suppose that the first and the second rows are equal, in particular it follows that $x_{i1}^{j} = x_{i2}^j$ for every $i=1,\ldots,n$ and for every $j=1,\ldots,d$. Since the action of $B$ just swaps rigidly the $n$ $N$-tuples of $(A^{N})^{n}$, any element $((y_{11},\ldots,y_{1N}),\cdots, (y_{n1},\ldots,y_{nN})) \tau$ of the subgroup generated by the $\gamma_j$'s satisfies $y_{11}=y_{12}$. This is a contradiction with assumption that the $\gamma_j$'s generate $G$.

Therefore $N\le \norm{A}^{nd}$ and applying logarithms on both sides of the inequality we have 
 \[
   d \ge \frac{1}{n\log\norm{A}} \log N = \frac{1}{n} \log_{\norm{A}} N > \frac{1}{n}\left( \mathrm{d}\left(A^N\right) - \d(A)-1\right),
 \]
where the last inequality holds by \cite[Lemma 2]{MR499355}.
\end{proof}
\end{lem}

The groups under study here are very different from the ones in \cite{MR2727305}. We cannot rely on the tree-like structure of iterated wreath products and the iterated exponentiation of permutation groups is not associative. This is the reason why we need to ask that the groups in the sequence $\mS$ have non-regular action. 
Before proving Theorem A we fix some notation. 

\phantom{a}

\noindent\bbb{Notation.} We will denote a $\wt{m}_k$-tuple in $\{1,\ldots,m_{k+1}\}^{\wt{m}_k}$ as $(i_1,\ldots,i_{\wt{m}_k})_{\wt{m}_k}$. This notation will be convenient in particular when we will have to deal with $\wt{m}_k$-tuples where all the coordinates are equal, for example $(1,\ldots,1)_{\wt{m}_k}$. We will denote an element of the group $S_{k+1}^{\wt{m}_k}$ as $(\sigma_1,\ldots,\sigma_{\wt{m}_k})_{\wt{m}_k}$.

 \begin{defn*}
 Let $(G,\{1,\ldots,m\})$ and $(H,\{1,\ldots,n\})$ be permutation groups. For $(i_1,\ldots,i_n),(j_1,\ldots,j_n) \in \{1,\ldots,m\}^n$ we say that $(i_1,\ldots,i_n)$ \bbb{precedes} $(j_1,\ldots,j_n)$, if and only if it exists $1\le l\le n$ such that $i_k=j_k$ $1\le k \le l-1$ and $i_{l} < j_{l}$. The relation ``precedes'' defines a total order on $\{1,\ldots,m\}^n$ that we will call \bbb{lexicographic order}. 
 \end{defn*}
 
 The following straightforward Lemma is one of the key tricks to prove Theorem A.
\begin{lem}\label{lem:trivdiag}
 Let $G \le \S(\ul{m})$ and $H \le \S(\ul{n})$ be permutation groups. Then the subgroup $H$ of the exponentiation $G\mbox{\textnormal{$\cwr$}} H$ acts trivially on the subset $$\{(i,\ldots,i) \vvert i\in \ul{m}\}.$$
\end{lem}
 
 Theorem A will now follow from an application of the next lemma.
 
\begin{lem}\label{lem:dgen}
 Let $\mathcal{S} = \{S_k \le \S(\ul{m}_k)\}_{k\in \ul{n}}$ be a sequence of transitive permutation groups and let $d$ be an integer. Suppose that $S_k$ is perfect, at most $d$-generated and for every $k \in \ul{n}$ there exist elements $i,j\in \ul{m}_k$ such that $\mathrm{St}_{S_k}(i) \neq \mathrm{St}_{S_k}(j)$. Then we exhibit a set of $d + \mathrm{d}(S_1)$ generators for $\wt{S}_n$. In particular $\wt{S}_n$ is generated by the generators of $S_1$ together with other $d$ special elements.
 \begin{proof}
 Without loss of generality we can suppose that $\mathrm{d}(S_k) = d$ for $k=2,\ldots,n$. Let $S_1 = \langle \a_1(1),\ldots,\a_{\mathrm{d}(S_1)}(1)\rangle$, $S_k = \langle \a_1(k),\ldots,\a_d(k)\rangle$, for $k=2,\ldots,n$ and order the elements of $\{1,\ldots,m_{k+1}\}^{\wt{m}_k}$ with lexicographic order. Without loss of generality we can suppose that for every $k\in \ul{n}$ we have 
\begin{equation}\label{eq:stab}
  \mathrm{St}_{S_k}(1)\neq \mathrm{St}_{S_k}(2).
\end{equation}

 We will now define $d$ elements of $\wt{S}_n$ that together with the generators of $S_1$ will generate $\wt{S}_n$. Define the elements $\beta_1,\ldots,\beta_d \in \wt{S}_n$ as
\[
     \b_j = (\a_j(n),e,\ldots,e)_{\wt{m}_{n-1}} \cdot (\a_j(n-1),e,\ldots,e)_{\wt{m}_{n-2}} \cdots (\a_j(2),e,\ldots,e)_{\wt{m}_1} 
\]
for $j=1,\ldots,d$. Note that the $\a_j(k)$'s are in the first place of the $\wt{m}_{k-1}$-tuples, which corresponds to the element $(1,\ldots,1)_{\wt{m}_{k-1}} \in \{1,\ldots,m_{k}\}^{\wt{m}_{k-1}}$.

Let $A \defeq \langle \a_1(1),\ldots,\a_{\mathrm{d}(S_1)}(1),\b_1,\ldots,\b_d \rangle \le \wt{S}_n$. We claim that $A = \wt{S}_n$. We will prove by induction on $k$ that $\wt{S}_k \le A$ for $k=1,\ldots,n$. Trivially $\wt{S}_1 = S_1 \le A$. Suppose by inductive hypothesis that $\wt{S}_{k} \le A$, we have to show that we can write any element of $\wt{S}_{k+1}$ as a product of the generators in $A$. Because $\wt{S}_{k+1} = S_{k+1}^{\wt{m}_k} \cdot \wt{S}_k$, it will suffice to show that $S_{k+1}^{\wt{m}_k} \le A$. 

By the transitivity of $S_k$ there is an element $t\in S_k$ such that $1^t = 2$ and by inductive hypothesis the element $\sigma \defeq (e,\ldots,e,t)_{\wt{m}_{k-1}} \in S_{k}^{\wt{m}_{k-1}}$ belongs to $A$. By Lemma \ref{lem:trivdiag} it follows that for $j=k,\ldots,n$
\begin{equation}\label{eq:1}
   (1,\ldots,1)_{\wt{m}_{j}}^{\sigma} = (1,\ldots,1)_{\wt{m}_{j}}
\end{equation}
 and from the definition of the exponentiation
\begin{equation}\label{eq:2}
    (1,\ldots,1)_{\wt{m}_{k-1}}^{\sigma} = (1^e,\ldots,1^e,1^t)_{\wt{m}_{k-1}} = (1,\ldots,1,2)_{\wt{m}_{k-1}}.
\end{equation}
 Moreover, since $\wt{S}_k \le A$, $\b_j' = (\a_j(n),e,\ldots,e)_{\wt{m}_{n-1}} \cdots (\a_j(k+1),e,\ldots,e)_{\wt{m}_k}$ belongs to $A$. Set $\gamma_j = [\sigma,\b_j']$, it is clear that $\gamma_j \in A$. 
Hence, by \eqref{eq:1}, $(\a_j(l),e,\ldots,e)_{\wt{m}_{l-1}}^{\sigma} = (\a_j(l),e,\ldots,e)_{\wt{m}_{l-1}}$ for $l=k+2,\ldots,n$ and, by \eqref{eq:2}, $(\a_j(k+1),e,\ldots,e)_{\wt{m}_{k}}^{\sigma} = (e,\a_j(k+1),e,\ldots,e)_{\wt{m}_{k}}$. Therefore $\gamma_j = (\a_j(k+1),\a_j(k+1)^{-1},e,\ldots,e)_{\wt{m}_k}$.

$S_k$ is transitive $\ul{m}_k$, thus, by the definition of exponentiation, $S_k^{\wt{m}_{k-1}}$ is also transitive on $\ul{m}_k^{\wt{m}_{k-1}}$. To conclude the proof it is sufficient to show that we can write any element of the form $(\lambda,e,\ldots,e)_{\wt{m}_k}$ in $S_{k+1}^{\wt{m}_k}$ as a word in the $\gamma_j$'s. We can then move $\lambda$ using the transitive action of $S_k^{\wt{m}_{k-1}}$. 

As $S_{k+1}$ is perfect it will be sufficient to prove that we can write any commutator $([\lambda_1,\lambda_2],e,\ldots,e)_{\wt{m}_k}$ as a word in the $\gamma_j$'s. By \eqref{eq:stab} there are $s \in S_k$ and $r\in \ul{m}_{\hspace{0.25mm}k}$, $r\neq 2$, such that $1^s = 1$ and $2^s = r$. By inductive hypothesis $\mu \defeq (e,\ldots,e,s)_{\wt{m}_{k-1}}$ belongs to $A$. Let $\lambda_1,\lambda_2 \in S_{k+1}$. Since the $\a_j(k+1)$'s generate $S_{k+1}$, there exist two $d$-variables words $w_1$ and $w_2$ such that $\lambda_1 = w_1(\a_1(k+1),\ldots,\a_d(k+1))$ and $\lambda_2 = w_2(\a_1(k+1),\ldots,\a_d(k+1))$. Thus, if we set $\delta_i = w_i(\a_1(k+1)^{-1},\ldots,\a_d(k+1)^{-1})$ for $i=1,2$, the elements $w_1(\gamma_1,\ldots,\gamma_d) = (\lambda_1,\delta_1,e,\ldots,e)_{\wt{m}_k}$ and $w_2(\gamma_1,\ldots,\gamma_d) = (\lambda_2,\delta_2,e,\ldots,e)_{\wt{m}_k}$ belong to $A$. The definition of $\mu$ and an easy calculation now yield 
\[
  \left[ (\lambda_1,\delta_1,e,\ldots,e)_{\wt{m}_k},(\lambda_2,\delta_2,e,\ldots,e)_{\wt{m}_k}^\mu \right] = ([\lambda_1,\lambda_2],e,\ldots,e)_{\wt{m}_k}
\]
Thus for every $\lambda \in S_{k+1}$ the $\wt{m}_k$-tuple $(\lambda,e,\ldots,e)_{\wt{m}_k}$ is in $A$. It follows that $S_{k+1}^{\wt{m}_k} \le A$ and $\wt{S}_{k+1} = S_{k+1}^{\wt{m}_k} \cdot \wt{S}_k \le A$. The result follows by induction. 
 \end{proof}
\end{lem}

\noindent \begin{namedproof}[Proof of Theorem A]
For every $n\in \mathbb{N}$, Lemma \ref{lem:dgen} gives us $d + \mathrm{d}(S_1)$ generators $\a_1(1),\ldots,\a_{\mathrm{d}(S_1)}(1),\b_1^{(n)},\ldots,\b_d^{(n)}$ of $\wt{S}_n$. For $n\in \mathbb{N}$, let $\pi_n$ be the inverse limit projection from $\varprojlim \wt{S}_k$ to $\wt{S}_n$. Let $a_1(1),\ldots,a_{\mathrm{d}(S_1)}(1),b_1,\ldots,b_d$ be the unique elements of $\varprojlim \wt{S}_k$ such that $\pi_n(a_{i}(1)) = \a_{i}(1)$ and  $\pi_n(b_j) = \beta_j^{(n)}$  for all $i\in \ul{\d(S_1)}$,$j\in \ul{d}$ and all $n \in \mathbb{N}$ and all $i=1,\ldots,d$. Then $a_1(1),\ldots,a_{\mathrm{d}(S_1)}(1),b_1,\ldots,b_d$ generate $\varprojlim \wt{S}_k$ by \cite[Proposition 4.1.1]{wilson:profinitegroups}. \qed
\end{namedproof}

Using \cite[Lemma 2]{segal:thefiniteimages} it is possible to improve the previous bound for 2-generated groups with the same method.

\begin{lem}\label{lem:3gen}
 Let $\mathcal{S} = \{S_k \le \S(\ul{m}_k)\}_{k\in \ul{n}}$ be a sequence of perfect 2-generated transitive permutation groups such that for every $k\in \ul{n}$ and every $i,j \in \ul{m}_k$ $\mathrm{St}_{S_k}(i) \neq \mathrm{St}_{S_k}(j)$. Then we exhibit a set of 3 generators for $\wt{S}_n$. In particular $\wt{S}_n$ is generated by the generators of $S_1$ together with another special element.
 \begin{proof}
 Let $S_k = \langle \a_1(k),\a_2(k) \rangle$. By \cite[Lemma 2]{segal:thefiniteimages}, for $k\in \ul{n}$, there exists $\sigma_k \in S_k$ and $1\le r_k\le m_k$ such that $r_k^{\sigma_k^2} \neq r_k$. Let
\[
  \beta = (\ldots,\a_1(n),\ldots,\a_2(n),\ldots)_{\wt{m}_{n-1}} \cdot \ldots \cdot (\ldots,\a_1(2),\ldots,\a_2(2),\ldots)_{\wt{m}_{1}} 
\]
where the element $\a_1(2)$ is in position $r_1^{\sigma_1}$, $\a_2(2)$ is in position $r_1$, $\a_1(k+1)$ is in position $(r_{k}^{\sigma_k},\ldots,r_{k}^{\sigma_k})_{\wt{m}_{k-1}}$, $\a_2(k+1)$ is in position $(r_{k},\ldots,r_{k})_{\wt{m}_{k-1}}$ for $k=2,\ldots,n-1$ and the identity in all the unspecified positions. Set $A \defeq \langle \a_1(1),\a_2(1),\b \rangle$ and proceed exactly as in the proof of Lemma \ref{lem:dgen}, with $\beta$ instead of $\beta_i$ and $(\sigma_k,\ldots,\sigma_k)_{\wt{m}_{k-1}}$ instead of $\sigma$, to show that $A=\wt{S}_n$.
\end{proof}
\end{lem}

Using Lemma \ref{lem:3gen} in place of Lemma \ref{lem:dgen} in the proof of Theorem A yields the following corollary.

\begin{cor}\label{cor:3gen}
 Let $\mathcal{S} = \{S_k \le \S(\ul{m}_k)\}_{k\in \mathbb{N}}$ be a sequence of perfect 2-generated transitive permutation groups. Suppose that for every $k\in \mathbb{N}$ and every $i,j \in \ul{m}_k$ we have $\mathrm{St}_{S_k}(i) \neq \mathrm{St}_{S_k}(j)$. Then we exhibit a set of 3 generators for $\varprojlim \wt{S}_n$.
\end{cor}

Finite generation of iterated wreath products of finite non-abelian simple transitive permutation groups was studied in \cite{quick:probabilisticgeneration} and the minimal number of generators is 2 in this case. We are convinced that in the case of perfect non-simple groups Lemma \ref{lem:3gen} is best possible but we do not have an explicit example to confirm this.

\begin{conj}
 There exists a sequence $\mS$ satisfying the hypotheses of Lemma \ref{lem:3gen} such that the infinitely iterated exponentiation of the groups in $\mS$ is not 2-generated.
\end{conj}

\section{Proof of Theorem B}\label{sec:B}
We now proceed to the proof of Theorem B. We will use the following.
\begin{thm*}(\cite[Theorem 3.1]{kaluzhnin})
 Let $n_1,n_2$ and $n_3$ be integers and let $A\le \S(\ul{n}_1), B\le \S(\ul{n}_2)$ and $C\le \S(\ul{n}_3)$ be permutation groups. Then $A \mbox{\textnormal{\cwr}} (B \wr C)$ and $(A\mbox{\textnormal{\cwr}} B) \mbox{\textnormal{\cwr}} C$ are isomorphic as permutation groups.
\end{thm*}
The next lemma is an application of \cite[Theorem 3.1]{kaluzhnin} and it will be used in the proof of Theorem B. 
\begin{lem}\label{lem:iso}
Let $\mS = \{S_k \le \S(\ul{m}_k)\}_{k\in \mathbb{N}}$ be a sequence of finite permutation groups. Let $\{k_n\}_{n\in \mathbb{N}}$ be an increasing sequence of integers and let $G_n$ be an iterated mixed wreath product of type $(\mS,\{k_n\}_{n\in \mathbb{N}})$. Define the permutation groups $\{\wh{S}_{k_n}^{(i)}\}$ for $n\in \mathbb{N}$ and $i\in \ul{k}_n\smallsetminus\ul{k}_{n-1}$ as follows: $\wh{S}_{k_n}^{(k_n)} = S_{k_n}$ and $\wh{S}_{k_n}^{(i)} = \wh{S}_{k_n}^{(i+1)} \mbox{\textnormal{$\cwr$}} S_{i}$. Set $H_i = \wh{S}_{k_i}^{(k_{i-1}+1)}$ for $i\in \ul{n}$.

Then $G_{k_n}$ is isomorphic (as a permutation group) to $\wt{H}_n$ for every $n\in\mathbb{N}$.
\begin{proof}
 The proof is by induction on $n$. If $n=1$ and $k_1=1$ the claim is trivial. If $k_1 >1$ repeated applications of \cite[Theorem 3.1]{kaluzhnin} yield the result. 
 
 Suppose that $G_{k_{n-1}}\cong \wt{H}_{n-1}$. By construction $G_{k_n}\cong S_{k_n}\cwr G_{k_n -1}$ and $G_{i} \cong S_i \wr G_{i-1}$ for $i\in \ul{k}_n \smallsetminus \ul{k}_{n-1}$. Therefore repeated applications of \cite[Theorem 3.1]{kaluzhnin} yield $G_{k_n} \cong (\wh{S}_{k_n}^{(i+1)} \cwr S_{i})\cwr G_{i-1}$ 
 for $i\in \ul{k}_n\smallsetminus\ul{k}_{n-1}$. 
Thus $G_{k_n} \cong  \wh{S}_{k_n}^{(k_{n-1}+1)}  \cwr G_{k_{n-1}}$ and by inductive hypothesis we conclude $G_{k_n} \cong H_n \cwr \wt{H}_{n-1} \cong \wt{H}_n$. 
The claim follows by induction.
\end{proof}
\end{lem}

\begin{lem}\label{lem:H}
 Let $A\le \S(\ul{m})$ and $B \le \S(\ul{n})$ be a permutation groups and set $G=A\mbox{\textnormal{\cwr}} B$. Then there exist $x,y \in \ul{m}^n$ such that $\mathrm{St}_{G}(x) \neq \mathrm{St}_G(y)$.
 \begin{proof}
  Consider the elements $x=(1,\ldots,1)_n$ and $y=(2,1,\ldots,1)_n$ in $\ul{m}^n$. Because $B$ is a permutation group there exists $b\in B$ such that $1^b = r \neq 1$, so $x^b = x$ and $y^b = (1,\ldots,1,2,1,\ldots,1)_n \neq y$. So $b$ is in the stabiliser of $x$ but not in the stabiliser of $y$.
 \end{proof}
\end{lem}

Finally we use Lemma \ref{lem:dgen}, Lemma \ref{lem:iso} and Lemma \ref{lem:H} to prove Theorem B.

\begin{namedproof}[Proof of Theorem B]
Let $G =\varprojlim G_n$ be an infinitely iterated mixed wreath product of type $(\mS,\{k_n\}_{n\in \mathbb{N}})$ and of stride at most $m$. We will use the same notation of Lemma \ref{lem:iso}. $G_{k_n}$ is isomorphic to $\wt{H}_n$ for every $n\in \mathbb{N}$, hence it is sufficient to show that the groups $\{H_n\}_{n\in \mathbb{N}}$ satisfy the hypotheses of Lemma \ref{lem:dgen}. It is clear that every $H_n$ is perfect and it can be generated by $md$ elements because it is an iterated wreath product of length at most $m$ made of $d$-generated groups. Moreover the action of $H_n$ is not regular by Lemma \ref{lem:H}. Applying Lemma \ref{lem:dgen} yields the result. \qed
\end{namedproof}

\begin{rmk}\label{rmk:weak}
 We can weaken the hypothesis of Theorem B in the following ways. $\{k_n\}_{n\in\mathbb{N}}$ is an increasing sequence of integers, $\mS$ is a sequence of finite perfect, at most $d$-generated, transitive permutation groups and: 
 \begin{itemize}
  \item for every $k_n$ such that $k_n = k_{n-1}+1$ there exist elements $i,j\in \ul{m}_{k_n}$ that have diffent stabilisers for the action of $S_{k_n}$. 
  \item $m\ge 2$.
 \end{itemize}
 The proof of Theorem B with these hypotheses remains the same.
\end{rmk}
 
\section{One application}\label{sec:special}
In this section we find explicitly two generators for infinitely iterated exponentiation of particular sequences $\mS$. We start with a lemma. 

\begin{lem}\label{lem:special}
 Let $\mathcal{S} = \{S_k \le \S(\ul{m}_k)\}_{k\in \ul{n}}$, be a sequence of finite 2-generated perfect transitive permutation groups. Suppose that for every $k \in \ul{n}$ there exist two generators $a_k,b_k$ of $S_k$ such that:
 \begin{itemize}
   \item $\mathrm{fix}(a_k)$ and $\mathrm{fix}(b_k)$ are non-empty,
   \item $(\norm{a_1},\norm{b_j}) = 1$ and $(\norm{b_1},\norm{a_j}) = 1$ for $j=2,\ldots,n$.
 \end{itemize}
 Then we exhibit two generators for $\wt{S}_n$.
 \begin{proof}
  Let $u_k \in \mathrm{fix}(a_k)$, $v_k \in \mathrm{fix}(b_k)$. 
  In the spirit of Lemma \ref{lem:dgen} define the following elements of $\ul{m}_i^{\wt{m}_{i-1}}$ 
  \[
    \underline{u}_i = (u_i,\ldots,u_i)_{\wt{m}_{i-1}} \mbox{\quad and \quad} \underline{v}_i = (v_i,\ldots,v_i)_{\wt{m}_{i-1}}
  \]
  for $i=2,\ldots,n-1$. By the transitivity of $S_k$ there is $\sigma \in S_k$ such that $u_{k}^{\sigma}= v_k$ and, by Lemma \ref{lem:trivdiag}, $\mu \defeq (\sigma,\ldots,\sigma)_{\wt{m}_{k}}$ is such that 
\begin{equation}\label{eq:a}
  \underline{u}_j^{\mu} = \underline{u}_j \mbox{\quad and \quad} \underline{v}_j^{\mu^{-1}} = \underline{v}_j 
\end{equation}
 for every $j\ge k+1$ and by definition of exponentiation we have
\begin{equation}\label{eq:b}
  \underline{u}_{\, k}^\mu = (u_{k}^{\sigma},\ldots,u_{k}^{\sigma})_{\wt{m}_{k}} = \underline{v}_{\; \!k}.
\end{equation}
For the rest of the proof we will write the position of an element in a tuple below the element itself. We claim that the elements
\begin{multline*}
\beta_1 = (e,\ldots,e,\underset{\underline{u}_{n-1}}{a_n},e\ldots,e)_{\wt{m}_{n-1}} \cdots (e,\ldots,e,\underset{\underline{u}_2}{a_3},e\ldots,e)_{\wt{m}_{2}}\cdot \\ \cdot (e,\ldots,e,\underset{v_1}{a_2},e\ldots,e)_{\wt{m}_{1}} b_1  
\end{multline*}
and 
\begin{multline*}
  \beta_2 = (e,\ldots,e,\underset{\underline{v}_{n-1}}{b_n},e\ldots,e)_{\wt{m}_{n-1}} \cdots (e,\ldots,e,\underset{\underline{v}_2}{b_3},e\ldots,e)_{\wt{m}_{2}}\cdot \\ \cdot (e,\ldots,e,\underset{u_1}{b_2},e\ldots,e)_{\wt{m}_{1}} a_1
\end{multline*}
generate the group $\wt{S}_n$. Let $A=\langle \b_1,\b_2\rangle$, we will prove by induction that $\wt{S}_k\le A$ for $k=1,\ldots,n$. It follows from Lemma \ref{lem:trivdiag} and the definition of $u_i$ and $v_i$ that $(e,\ldots,e,{a_{i}},e\ldots,e)_{\wt{m}_{i}}$ commutes with $(e,\ldots,e,{a_{j}},e\ldots,e)_{\wt{m}_{j}}$ for $i\neq j$. Set $p=\prod_{i=2}^n \norm{a_i}$ and $q=\prod_{i=2}^n \norm{b_i}$, then $\b_1^p = b_1^p$ and $\b_2^q = a_1^q$, so $S_1 \le A$.

Suppose by inductive hypothesis that $\wt{S}_k\le A$, our goal is to write any element of $S_{k+1}^{\wt{m}_k}$ as a word in $\b_1$, $\b_2$. Clearly the elements
\[
   \beta_1' = (e,\ldots,e,\underset{\underline{u}_{n-1}}{a_n},e\ldots,e)_{\wt{m}_{n-1}} \cdots  (e,\ldots,e,\underset{\underline{u}_{k}}{a_{k+1}},e\ldots,e)_{\wt{m}_{k}}
\]
and
\[
  \beta_2' = (e,\ldots,e,\underset{\underline{v}_{n-1}}{b_n},e\ldots,e)_{\wt{m}_{n-1}} \cdots  (e,\ldots,e,\underset{\underline{v}_{k}}{b_{k+1}},e\ldots,e)_{\wt{m}_{k}} 
\]
belong to $A$.

Let us now consider the commutators $\gamma_i \defeq [\mu_i,\b_i']$ for $i=1,2$. Following exactly the steps of Lemma \ref{lem:dgen} we can use \eqref{eq:a}, \eqref{eq:b} (instead of \eqref{eq:1} and \eqref{eq:2}) and Lemma \ref{lem:trivdiag} to show $S_{k+1}^{\wt{m}_k} \le A$. Therefore $\wt{S}_{k+1} = S_{k+1}^{\wt{m}_k} \cdot \wt{S}_k \le A$. The result follows by induction.
\end{proof}
\end{lem}

Since none of the groups $\wt{S}_k$ is cyclic, the proof that $\varprojlim \wt{S}_k$ is topologically 2-generated is now the same as the proof of Theorem A using Lemma \ref{lem:special} instead of Lemma \ref{lem:dgen}. We proved the following.

\begin{cor}\label{cor:special}
 Let $\mathcal{S} = \{S_k \le \S(\ul{m}_k)\}_{k\in \mathbb{N}}$, be a sequence of finite 2-generated perfect transitive permutation groups. Suppose that for every $k \in \ul{n}$ there exist two generators $a_k,b_k$ of $S_k$ such that:
 \begin{itemize}
   \item $\mathrm{fix}(a_k)$ and $\mathrm{fix}(b_k)$ are non-empty,
   \item $(\norm{a_1},\norm{b_j}) = 1$ and $(\norm{b_1},\norm{a_j}) = 1$ for $j=2,\ldots,n$.
 \end{itemize}
 Then the infinitely iterated exponentiation $\varprojlim \wt{S}_k$ is topologically 2-generated and we produce explicitly two generators for the group.
\end{cor}

\begin{rmk}\label{rmk:most}
 Using the Classification of Finite Simple Groups it was shown that all finite non-abelian simple groups besides $Sp_4(2^n)$, $Sp_4(3^n)$, $\phantom{}^2B_2(2^{2n+1})$ and possibly other finitely many exceptions can be generated by an involution and an element of order 3, see \cite[Corollary 1.3]{MR1688493} for further references on $(2,3)$-generation of finite non-abelian simple groups. A sequence of $(2,3)$-generated non-abelian simple groups satisfies the hypotheses of Lemma \ref{lem:special}. 
\end{rmk} 

\section*{Acknowledgements}
This work was carried out as part of my PhD studies at Royal Holloway University of London under the supervision of Dr.~Yiftach Barnea, I would like to thank him for introducing me to the problem and his continuous help and support. I also wish to thank Eugenio Giannelli and Benjamin Klopsch for the most useful discussions, in particular I am grateful to the latter for the suggestion of Lemma \ref{lem:lowerbound}.

\bibliographystyle{plain}
\def\cprime{$'$}

\vspace{1cm}
 
\noindent MATTEO VANNACCI

\noindent Royal Holloway, University of London\\
Egham, Surrey, TW20 0EX\\
United Kingdom

\noindent\texttt{vannacci.m@gmail.com}
\end{document}